\documentclass[review]{elsarticle}

\usepackage{lineno,hyperref}
\modulolinenumbers[5]

\usepackage{amsmath}
\usepackage{amssymb}
\usepackage{graphicx}
\usepackage{subfig}
\usepackage{float}
\usepackage{hyperref}
\usepackage{amssymb,amsmath,amsthm}
\usepackage{color}

\newtheorem{thm}{Theorem}
\newtheorem{lem}[thm]{Lemma}

\newtheorem{prob}[thm]{Problem}

\newdefinition{rem}{Remark}
\newdefinition{exa}{Example}
\newproof{pf}{Proof}
\newproof{pot}{Proof of Theorem \ref{thm2}}
\newtheorem{algorithm}{Algorithm}[section]{\bf}{\rm}
\def\trace{\mathop{\mathrm{trace}}}

\def\grad{\mathop{\mathrm{grad}}}
\def\argmin{\mathop{\mathrm{argmin}}}
\def\vec{\mathop{\mathrm{vec}}}
\journal{Journal of Computational and Applied Mathematics}

\usepackage{algorithm}
\usepackage[noend]{algpseudocode}
\usepackage{breakcites}
\algrenewcommand{\algorithmiccomment}[1]{{\hfill\color{black}$\blacktriangleright$ #1}}

\begin{document}
	
	\begin{frontmatter}	
		
\title{Solving the Discrete Euler-Arnold Equations \\ for the Generalized Rigid Body Motion}

\author{Jo\~ao R. Cardoso\corref{mycorrespondingauthor}}
\address{Coimbra Polytechnic--ISEC, and \\ Center for Mathematics, University of Coimbra, Portugal.}
\ead{jocar@isec.pt}

\cortext[mycorrespondingauthor]{Corresponding author: Jo\~ao R. Cardoso}

\author{Pedro Miraldo\corref{}}
\address{Institute for Systems and Robotics (LARSyS), \\ Instituto Superior T\'ecnico, University of Lisbon, Portugal.}
\ead{pedro.miraldo@tecnico.ulisboa.pt}

\begin{abstract}
We propose three iterative methods for solving the Moser-Veselov equation, which arises in the discretization of the Euler-Arnold differential equations governing the motion of a generalized rigid body. We start by formulating the problem as an optimization problem with orthogonal constraints and proving that the objective function is convex. Then, using techniques from optimization on Riemannian manifolds, the three feasible algorithms are designed. The first one splits the orthogonal constraints using the Bregman method, whereas the other two methods are of the steepest-descent type. The second method uses the Cayley-transform to preserve the constraints and a Barzilai-Borwein step size, while the third one involves geodesics, with the step size computed by Armijo's rule. Finally, a set of numerical experiments are carried out to compare the performance of the proposed algorithms, suggesting that the first algorithm has the best performance in terms of accuracy and number of iterations. An essential advantage of these iterative methods is that they work even when the conditions for applicability of the direct methods available in the literature are not satisfied.
\end{abstract}


\end{frontmatter}


\section{INTRODUCTION}

In \cite{Moser91}, Moser and Veselov proposed the following equations to discretize the classical Euler-Arnold differential equations for the motion of a generalized rigid body:
\begin{equation}
\begin{aligned}
\label{eq:euler_arnold}
    M_{k+1} & = \omega_k M_k\omega_k^T \\
    M_k &= \omega_k^T J-J\omega_k,
\end{aligned}
\end{equation}
where $M_k$ is the angular momentum with respect to the body (here represented by a skew-symmetric matrix), $J$ is the inertia matrix (symmetric positive definite), and $\omega_k$ (orthogonal matrix) is the angular velocity. Rigid body equations arise in several applications, e.g., celestial mechanics, molecular dynamics, mechanical robotics, and flight control, where they are used in particular to understand the body--body interactions of particles like planets, atoms, and molecules. See, for instance, \cite{Bloch15,Dang20,Lee07,Nordkvist10,Kalabic14,Kalabic15} and the references therein.

The main challenge of solving (\ref{eq:euler_arnold}) is to find an orthogonal matrix $\omega_k$ in the second equation, by assuming that $J$ and $M_k$ are given. Mathematically, the problem consists of finding an orthogonal matrix ${X}$ (for convenience, here we use $X$ instead of $\omega$) such that 
\begin{equation}
    \label{eq:moser_veselov}
    {X}{J} - {J}{X}^T = {M},
\end{equation}
where ${J}$ is a given symmetric positive definite matrix, and ${M}$ is a known skew-symmetric matrix. All the matrices involved are square of order $n$. The matrix equation (\ref{eq:moser_veselov}) is known as the {\it Moser-Veselov} equation and was firstly investigated in \cite{Moser91}, where the authors based their developments on factorizations of certain matrix polynomials. A different approach, but computationally more efficient, was provided later in \cite{Cardoso03}, where the authors noted that (\ref{eq:moser_veselov}) can be connected with a certain algebraic Riccati equation and, in turn, with the {\it Hamiltonian} matrix
\begin{equation}\label{eq:hamiltonian}
{\cal H}=\left[
\begin{array}{cc}
M/2&I\\M^2/4+J^2&M/2
\end{array}
\right].
\end{equation}
We now revisit some results stated in \cite{Cardoso03}, concerning the existence and uniqueness of solutions of \eqref{eq:moser_veselov}.

\begin{thm}\label{thm-existence} 
For the matrix equation \eqref{eq:moser_veselov}:
\begin{enumerate}
	\item There exists a solution $X\in \mathcal{SO}(n)$ (the special orthogonal or rotation group of order $n$)  if and only if the size of the Jordan blocks associated to the pure imaginary eigenvalues of ${\cal H}$ (if any) is even;
	\item \eqref{eq:moser_veselov} has a unique solution $X\in \mathcal{SO}(n)$ if and only if the spectrum of ${\cal H}$ is pure imaginary and the size of the Jordan blocks associated to each (nonzero) eigenvalue is even.
\end{enumerate}
\end{thm}

From Theorem~\ref{thm-existence}, we can see that the Moser-Veselov equation may have no solution in $\mathcal{SO}(n)$ if the associated Hamiltonian matrix ${\cal H}$ has any pure imaginary eigenvalue with a Jordan block of odd size. It is also known that the existence of purely imaginary eigenvalues in ${\cal H}$ causes significant difficulties in solving \eqref{eq:moser_veselov}. To avoid those situations, as far as we know, in all the existing algorithms for solving the equation, it is assumed a priori that ${\cal H}$ does not admit any pure imaginary eigenvalue (see, \cite[Sec. 1.2]{Mclachlan05}, \cite[Sec. 1.4]{Moser91}, and \cite[Sec. G]{Nordkvist10}). Moreover, the algorithms based on solving the associated algebraic Riccati equation require the strong condition that the matrix $M^2/4+J^2$ must be symmetric positive definite. These issues have motivated us to investigate methods whose applicability does not require those restrictive conditions. As we can see later in Sec.~\ref{experiments}, the three proposed optimization algorithms produce special orthogonal solutions, even when $M^2/4+J^2$ is not symmetric positive definite. Those iterative algorithms may also be used in problems where ${\mathcal H}$ has purely imaginary eigenvalues associated with Jordan blocks of even size (check Theorem \ref{thm-existence}) but, as will be illustrated later in Sec.~\ref{experiments}, the convergence may slow down.   

\begin{prob}\label{main-problem}
Let $\|.\|_F$ denote the Frobenius norm, i.e., $\|A\|_F:=\sqrt{\trace(A^TA)}$. The problem of finding a special orthogonal solution $X$ in \eqref{eq:moser_veselov} can be formulated as an optimization problem in the following way: 
\begin{equation}\label{problem1}
\min_{X\in \mathcal{SO}(n)} \left\|XJ - JX^T - M\right\|^2_F.
\end{equation}
\end{prob}

In Problem~\ref{main-problem}, we have chosen the Frobenius norm because its definition in terms of the trace of a matrix allows us to access the derivatives of the objective function easily, making it more suitable to handle optimization problems than other norms, like, for instance, spectral or infinity norms.

The literature on numerical methods for solving non-linear constrained problems, like \eqref{problem1}, is large; see for instance \cite{Nocedal06,Luenberger15,Polak97,Cardoso2017,Campos2019}. However, due to the complicated expression of the objective function and the large number of constraints arising from the conditions $X^TX=I$ and $\det(X)=1$, some care must be taken with the choice of the methods. 

Techniques from Riemannian geometry for solving optimization problems with orthogonal constraints have attracted the interest of many researchers in the last decades; see \cite{Absil07,Edelman99}, and the references therein. 
An essential feature of those techniques is that they allow the transformation of a constrained optimization problem into an unconstrained one. Moreover, since the set of orthogonal matrices is a manifold and provided that the objective function satisfies some smoothness requirements, we can make available tools such as Euclidean gradients, Riemannian gradients, retractions, and geodesics. 

The three methods presented in this work evolve on the orthogonal manifold and belong to the family of line search methods on manifolds described in \cite[Ch. 4]{Absil07}. 
They are iterative and feasible (or constraint-preserving), in the sense that, starting with a matrix $X_0\in \mathcal{SO}(n)$, all the iterates $X_k$ also stay in $\mathcal{SO}(n)$. 

The major contributions of this work are:
\begin{itemize}
\item the development of three effective algorithms for solving Problem \ref{main-problem} and in turn the matrix equation \eqref{eq:moser_veselov} (i.e., the Moser-Veselov equation), that do not require the condition $M^2/4+J^2>0$;
\item a detailed discussion on the properties of the optimization problems, including convexity issues; 
\item a novel approach for solving the unconstrained optimization problems arising in each iteration of the Bregman splitting algorithm and a discussion about the reasons that led the MATLAB's \texttt{fminunc} function to give unsatisfactory  results in some circumstances;
\item a novel relative residual capable to infer about the quality of the computed solution of the Moser-Veselov equation;
\item careful modifications on existing algorithms for solving optimization problems with orthogonal constraints, to make them suitable for our particular problems.
\end{itemize}

In the next section, we derived a workable expression to the objective function of Problem~\ref{main-problem}. Sec.~\ref{sec:algs} presents the three algorithms proposed in this paper. In Sec.~\ref{sec:numIssues}, numerical issues of the algorithms are discussed, and, in Sec.~\ref{experiments}, a selection of experiments are carried out to illustrate the performance of the algorithms. In Sec.~\ref{sec:conclusions}, some conclusions will be drawn.
\section{Rewriting the Objective Function}
Let us denote by $F(X):= \left\|XJ - JX^T - M\right\|^2_F$ the objective function arising in \eqref{problem1}. Using the properties of the trace of a matrix and attending that $J^T=J$ and $M^T=-M$, we have (detailed calculation is omitted):
\begin{align}
	F(X) &= 
	\trace\left( (XJ - JX^T - M)^T(XJ - JX^T - M)\right) \nonumber \\
	& =2\trace\left(JX^TXJ\right)-2\trace(XJXJ) + 
	 4\trace(MXJ) -\trace(M^2).
	\label{F}
\end{align}
Now, if we take into account the orthogonality of $X$, that is, $X^TX=XX^T=I$, then $F(X)$ can be simplified to (a different notation is used):
\begin{equation}\label{F-tilde}
	\widetilde{F}(X) =-2\trace\left((JX)^2\right) + 4\trace(XJM) + 2\trace(J^2) - \trace(M^2),
\end{equation}
which is the restriction of $F(X)$ to the orthogonal group $\mathcal{O}(n)$, that is, $F(X)=\widetilde{F}(X)$, for any $X\in \mathcal{O}(n)$, but, in general, $F(X)\neq\widetilde{F}(X)$, if $X\notin \mathcal{O}(n)$. Hence, the problem (\ref{problem1}) may be simplified to:
\begin{equation}\label{problem}
\min_{X\in \mathcal{SO}(n)} \widetilde{F}(X)=-2\trace\left((JX)^2\right)+4\trace(XJM)+\alpha,
\end{equation}
where $\alpha:= 2\trace(J^2)-\trace(M^2)$.

If $x_{ij}$ denotes the entry $(i,j)$ of the matrix $X$, then both $F(X)$ defined in (\ref{F}) and $\widetilde{F}(X)$ in (\ref{F-tilde}) are differentiable functions in $\mathbb{R}^{n^2}$, because they are quadratic polynomials in the $n^2$ variables $x_{ij}$. In the following lemma we show that $F(X)$ is convex in the set of all $n\times n$ matrices with real entries $\mathbb{R}^{n\times n}$, that is: \begin{equation}\label{def-convex}
F\left(tX_1+(1-t)X_2\right) \leq tF(X_1)+(1-t)F(X_2),
\end{equation}
for any $t\in [0,1]$ and $X_1,X_2\in \mathbb{R}^{n\times n}$.

\begin{lem}\label{convex}
	The function $F(X)$ given in (\ref{F}) is convex in  $\mathbb{R}^{n\times n}$.
\end{lem}

\begin{proof}
Let us denote $f(X):=\|XJ - JX^T - M\|_F$. Attending that the Frobenius norm satisfies the triangle inequality, we have
\begin{align}
		f\left(tX_1+(1-t)X_2\right) &= \left\|\left(tX_1+(1-t)X_2\right)J-J\left(tX_1+(1-t)X_2\right)^T-M\right\|_F \nonumber \\
		&= \left\|t(X_1J-JX_1^T)+(1-t)(X_2J-JX_2^T)-M\right\|_F \nonumber \\
		&= \left\|t(X_1J-JX_1^T)+(1-t)(X_2J-JX_2^T)-(t+1-t)M\right\|_F \nonumber \\
		&= \left\|t(X_1J-JX_1^T-M)+(1-t)(X_2J-JX_2^T-M)\right\|_F \nonumber \\
		&\leq t\left\|X_1J-JX_1^T-M\right\|_F+(1-t)\left\|X_2J-JX_2^T-M\right\|_F \nonumber  \\
		&= tf(X_1)+(1-t)f(X_2),
\end{align}
for all $t\in [0,1]$ and $X_1,X_2\in \mathbb{R}^{n\times n}$. Consider the scalar function $g(y)=y^2$. Since $f(X)\geq 0$, for any $X\in \mathbb{R}^{n\times n}$, and $g$ is non-decreasing in $[0,+\infty[$, we conclude that  $F(X)=g(f(X))$ is convex.
\end{proof}

Similarly, we could show that $\widetilde{F}(X)$ is convex in $\mathbb{R}^{n\times n}$. Note, however, that the constraints of the optimization problem (\ref{problem}) are non-convex, that is, for $P,\,Q\in \mathcal{SO}(n)$ and $t\in [0,1]$, in general $tP+(1-t)Q\notin \mathcal{SO}(n)$, which makes the problem much more difficult. 

Now, we use the rules for the derivatives of the trace function (see, for instance, \cite[Ch. 10]{Lutkepohl97}) to obtain the expressions of the Euclidean gradients (derivatives with respect to $X$) of those functions:
\begin{align}
\nabla F(X)&= 4XJ^2-4JX^TJ-4MJ\label{grad-F}\\
\nabla\widetilde{F}(X)&= -4JX^TX-4MJ.\label{grad-F-tilde}
\end{align}
The Riemannian gradients in the orthogonal manifold can be defined by:
\begin{eqnarray}
\grad F(X)&=& \nabla F(X)\, X^T - X\, \nabla F(X)^T\label{R-grad-F}\\
\grad\widetilde{F}(X)&=& \nabla\widetilde{F}(X)\, X^T-X\,\nabla\widetilde{F}(X)^T.\label{R-grad-F-tilde}
\end{eqnarray}
We recall that these Riemannian gradients belong to the orthogonal group's tangent space, which is the set of skew-symmetric matrices (see, for instance, \cite{Edelman99}).  

In the next section, we propose three algorithms to solve \eqref{main-problem}.

\section{Algorithms}\label{sec:algs}

There is a vast literature on methods for optimizing functions on the orthogonal group or, more generally, on the Stiefel manifold (e.g., \cite{Abrudan08,Edelman99,Gao18,Jiang15,Lai14,Manton02,Wen13,Zhu17}). Among those methods, we have selected three state-of-the-art ones that we believe to be well suited for our specific objective function. 
The first algorithm (Section~\ref{sec-alg1}) splits the orthogonal constraints in a Bregman's style \cite{Bregman67,Osher05} and is based on the SOC algorithm proposed by Lai and Osher in \cite[Alg. 2]{Lai14}. 
Our second algorithm (Section~\ref{sec-alg2}) is inspired in the feasible method developed by Wen and Yin in \cite[Alg. 2]{Wen13}, which uses a retraction --- the Cayley transform --- instead of geodesics. It is a line search method and, to find the appropriate step size, the method of Barzilai-Borwein (BB) \cite{Barzilai88} is used. The third algorithm involves line search techniques, namely the Armijo's rule, and is based on a proposal by Abrudan {\it et al.} in \cite[Table II]{Abrudan08}. It is of steepest descent type and involves geodesics, more specifically matrix exponentials. 

\subsection{Algorithm Based on Bregman Splitting}\label{sec-alg1}

The SOC Algorithm \cite[Alg. 2]{Lai14} applied to Problem~\ref{main-problem} is summarized in the following steps: 
\begin{enumerate}
	\item Choose a positive scalar $r$ and a starting matrix $X_0$. Set $P_0=X_0$ and $B_0=0$;
	\item While ``not converge'' do
	\begin{enumerate}
		\item $\displaystyle X_k=\argmin_X\ F(X)+\frac{r}{2}\|X-P_{k-1}+B_{k-1}\|_F^2,$ where $F(X)$ is defined by (\ref{F});
		\item $Y_k\leftarrow X_k+B_{k-1}$;
		\item Compute the singular value factorization: $Y_k=UDV^T$;
		\item $P_k\leftarrow UV^T$;
		\item $B_k\leftarrow B_{k-1}+X_k-P_k$.
	\end{enumerate}
\end{enumerate}

A drawback of this algorithm is the requirement of solving the (unconstrained) optimization problem in Step 2(a). Let us write down the corresponding objective function in terms of the trace.
Set $G_k(X):=\frac{r}{2}\|X-P_{k-1}+B_{k-1}\|_F^2$ and $C_k=-P_{k-1}+B_{k-1}$. After a few calculations, we have
\begin{equation}\label{G}
G_k(X)=\frac{r}{2}\left( \trace(X^TX)+2\trace(C_k^TX)+\trace(C_k^TC_k)\right).
\end{equation}
Hence, the objective function in Step 2(a) is ${\cal F}_k(X):=F(X)+G_k(X)$, where $F(X)$ denotes the function given in \eqref{F}, and the associated unconstrained optimization problem may be formulated as 
\begin{multline}
\min_X  {\cal F}_k(X)= \min_X \ 2\trace\left(JX^TXJ\right) - 2\trace(XJXJ) + 4\trace(MXJ) - \\ \trace(M^2) + 
\frac{r}{2}\left( \trace(X^TX) +  2\trace(C_k^TX)+\trace(C_k^TC_k)\right),\label{step2a}
\end{multline}
where $C_k,J,M$ are given square matrices of order $n$, $J$ and $M$ are, respectively, symmetric positive definite and skew-symmetric, and $r$ is a positive parameter.
By a similar argument to the one used in Lemma~\ref{convex}, we can show that, for each $k$, $G_k(X)$ in \eqref{G} is convex, and the same is valid to the objective function
${\cal F}_k(X)$. This is a very useful property because it guarantees that local minima are global as well.

Since each iteration $k$ of the above SOC algorithm requires the solution of a convex unconstrained optimization problem of the form \eqref{step2a}, whose objective function changes for each $k$, according to the entries of the matrix $C_k$, a possible approach to solve each one is to use the MATLAB's function \texttt{fminunc}, which is based on quasi-Newton and trust region methods. More precisely,  \texttt{fminunc} is based, by default, on the Broyden, Fletcher, Goldfarb, and Shanno quasi-Newton method, which is also know as the BFGS method (check \cite[Sec. 8.1]{Nocedal06} and the references therein). However, if the gradient $\nabla {\cal F}(X_k)$ is provided, then \texttt{fminunc} switches to a trust region method based on the proposals of Coleman and Li \cite{Coleman94,Coleman96}. 

We have solved many unconstrained problems of the form \eqref{step2a} using \texttt{fminunc} but, despite the convexity, the results were not so good as expected, either in terms of speed of convergence or in terms of accuracy. We learnt from our experiments that \texttt{fminunc} (without gradient and with the default tolerance of $10^{-6}$) is reliable only for very small size problems (say, $n\leq 5$). As $n$ increases, we observed that in some iterations $k$ of the SOC algorithm, \texttt{fminunc} was unable to minimize ${\cal F}(X_k)$. It displayed warnings like ``local minimum possible'' or ``solver stopped prematurely''. Recall that finding the minimum of \eqref{step2a} in all the iterations is necessary to guarantee the convergence of the SOC algorithm. For a given fixed $k$, let $X_k^{(i)}$ denotes the $i$-th iteration generated by \texttt{fminunc} when applied to the minimization  ${\cal F}(X_k)$. An interesting fact we have observed was that, as $i$ increased, ${\cal F}(X_k^{(i)})$ decreased quite fast to values lower than $10^{-6}$, while the components of the gradient vector $\nabla {\cal F}(X_k^{(i)})$ decreased in a slow fashion towards zero. This implies that \texttt{fminunc} involves a large number of iterations to guarantee that the norm of the gradient is lower than a fixed tolerance and hence that $X_k^{(i)}$ satisfies (up to that tolerance) the first-order necessary conditions. Recall that, if $X_\ast$ is a local (or global) minimizer of and ${\cal F}$ (which is continuously differentiable), we must have  $\nabla{\cal F}(X_\ast)=0$.
   
However, if the gradient function $\nabla{\cal F}$ (check \eqref{grad-cal-F}) is provided to  \texttt{fminunc}  and the tolerance is set to $10^{-5}$, we see that its performance improves and, 
as shown later in Section~\ref{experiments}, for equations of small size (say, $n\leq 15$), the usage of \texttt{fminunc} can be viewed as a possible approach to make the SOC algorithm effective.

To deal with smaller tolerances and equations involving matrices with larger size, we propose a different approach, which is described next. As we will see later, this approach seems to be very promising, either in terms of accuracy and computational cost. 

We start by finding the zeros of the gradient of ${\cal F}(X)$ (here, to simplify the notation, we omit the subscript $k$). We note that the expression of ${\cal F}(X)$ is non-linear and involves $n^2$ variables.

A few calculations lead to the following expression to the gradient of the objective function ${\cal F}(X)$:
\begin{equation}\label{grad-cal-F}
\nabla{\cal F}(X)=4XJ^2-4JX^TJ-4MJ+r(X+C).
\end{equation}
We know that local minima of ${\cal F}(X)$ are among the zeros of its Euclidean gradient. Since ${\cal F}$ is convex, those local minima (if any) will be global as well. Therefore, we need to investigate the solutions of the matrix equation  $\nabla{\cal F}(X)=0$, which is equivalent to 
\begin{equation}\label{matrix-eq-1}
X(4J^2+rI)-4JX^TJ=4MJ-rC.
\end{equation}
An easy way of solving \eqref{matrix-eq-1} is achieved by performing vectorization. Let $\vec(.)$ stand for the operator that stacks the columns of a matrix $n\times n$ into a long vector of size $n^2\times 1$, and let $\otimes$ denote the Kronecker product. It is well-known that $\vec(AYB)=(B^T\otimes A)\vec(Y)$ and that
\begin{equation}\label{permut}
\vec(A^T)=\Pi\vec(A),
\end{equation}
where $\Pi$ is the commutation (or permutation) matrix of order $n^2\times n^2$ (check \cite[Ch. 7, Sec. 9.2]{Lutkepohl97}). Applying the $\vec$ operator to (\ref{matrix-eq-1}) yields
\begin{align}
	\vec\left(X(4J^2+rI)\right)-4\vec(JX^TJ)&=\vec(4MJ-rC)\iff\nonumber\\
	 \left((4J^2+rI)\otimes I\right)\vec(X)-4(J\otimes J)\vec(X^T)&=\vec(4MJ-rC)\iff\nonumber\\
	\left((4J^2+rI)\otimes I\right)\vec(X)-4(J\otimes J)\Pi\vec(X)&=\vec(4MJ-rC)\iff\nonumber\\
	\left[(4J^2+rI)\otimes I-4(J\otimes J)\Pi\right]\vec(X)&=\vec(4MJ-rC).\label{eq-kron}
\end{align}
Note that (\ref{eq-kron}) corresponds to a linear system of type
$$A{\mathbf x}={\mathbf b},$$
where $A:=(4J^2+rI)\otimes I\,-\,4(J\otimes J)\Pi$ is an $n^2\times n^2$ matrix, ${\mathbf x}:=\vec(X)$, and ${\mathbf b}:=\vec(4MJ-rC)$. The vectorization approach is very useful to understand the theory of the matrix equation (\ref{matrix-eq-1}). If $A$ is non-singular, a unique solution to (\ref{matrix-eq-1}) is guaranteed. Otherwise, if $A$ is singular, (\ref{matrix-eq-1}) may have infinitely many solutions or no solutions. In the numerical examples we have considered (including the ones to be shown in Section \ref{experiments}), the matrix $A$ encountered was always non-singular.      

It turns out, however, that solving \eqref{matrix-eq-1} through the linear system $A{\mathbf x}={\mathbf b}$ would require $O(n^6)$ operations, which is prohibitive, especially if $n$ large. This has motivated us to look for less expensive methods for solving \eqref{matrix-eq-1}.

Since $J$ is assumed to be non-singular, a right multiplication of the matrix equation \eqref{matrix-eq-1} by $J^{-1}$ yields
 \begin{equation}\label{matrix-eq-2}
 X(4J+rJ^{-1})-4JX^T=4M-rCJ^{-1}.
 \end{equation}
Setting $Y:=X^T$, $A_1:=-4J$, $A_2:=4J+rJ^{-1}$, and $A_3:=4M-rCJ^{-1}$, \eqref{matrix-eq-2} can be rewritten in the form 
 \begin{equation}\label{matrix-eq-3}
A_1Y+Y^TA_2=A_3,
\end{equation}
which is a Sylvester-type equation. An effective method to solve \eqref{matrix-eq-3} may be found in \cite[Alg. 3.1]{Teran11}, which is based on the so-called QZ decomposition. Do not confuse \eqref{matrix-eq-3} with the more easily to handle classical Sylvester matrix equation $A_1X+XA_2=A_3$. Although the QZ decomposition is quite expensive, it can be performed in $O(n^3)$ operations. Hence, it is possible to solve \eqref{matrix-eq-1} and the Step 2(a) of the SOC algorithm efficiently by $O(n^3)$ operations, which is the typical cost for problems involving matrix-matrix products.

A drawback we have noticed when performing experiments with a direct application of the algorithm of \cite{Lai14}, to solve the Moser-Veselov equation, was the poor feasibility of the approximation obtained when compared to the other algorithms to be addressed in the next sections. That is, denoting by $\widetilde{X}$ the approximation obtained for the solution, we have observed that the value $\| \widetilde{X}^T\widetilde{X} -I\|_F$ was not close enough to zero. To overcome this issue, we project $X_k$ onto the orthogonal group by computing its singular value decomposition (SVD). Provided that \eqref{matrix-eq-3} has a unique solution $Y_\ast$, that means that $X_\ast=Y_\ast^T$ is the unique minimizer of the objective function in \eqref{step2a}.

Algorithm \ref{alg:bregman} summarizes the main steps of the Bregman splitting algorithm to solve the Moser-Veselov equation, \eqref{eq:moser_veselov}. The next subsection presents two alternative approaches for solving this problem.

\begin{algorithm}[t]
	\caption{Bregman splitting algorithm to solve the Moser-Veselov equation \eqref{eq:moser_veselov}.}
	\label{alg:bregman}
	{
		\begin{algorithmic}[1]
			\State Choose a positive scalar $r$ and a starting matrix $X_0\in \mathcal{SO}(n)$;
			\State $P_0\gets X_0$; $B_0\gets 0$;
			\State $J_1\gets J_1^{-1}$; $A_1\gets -4J$; $A_2\gets 4J+rJ_1$; 
			\State $k\gets 1$; 
			\While  {``not converge''} 
			\State $A_3\gets 4M-r(B_{k-1}-P_{k-1})J_1$;
			\State Solve $A_1X_k+X_k^TA_2=A_3$ for $X_k$, using \cite[Algs. 3.1]{Teran11};
			\State $X_k\gets X_k^T$;
			\State $Y_k\gets X_k+B_{k-1}$;
			\State Compute the SVD $Y_k=UDV^T$ and set $P_k\gets UV^T$;
			\State Compute the SVD $X_k=UDV^T$ and set $X_k\leftarrow UV^T$;
			\State $B_k\gets B_{k-1}+X_k-P_k$;
			\State $k\gets k+1$.
			\EndWhile
		\end{algorithmic}
	}
\end{algorithm}

\subsection{Two Steepest Descent-Type Algorithms}\label{sec-alg2}

A successful method to solve optimization problems with orthogonal constraints is the feasible iterative method developed in \cite[Alg. 2]{Wen13}. At each iteration, the skew-symmetric Riemannian gradient \eqref{R-grad-F-tilde} is multiplied by a suitable positive number $\tau$ (step-size), and transformed into an orthogonal matrix by means of the Cayley transformation
\begin{equation}\label{cayley}
Y(\tau)=\left(I+\frac{\tau}{2}W\right)^{-1} \left(I-\frac{\tau}{2}W\right)X,
\end{equation}
where $W$ is the Riemannian gradient and $X$ is an orthogonal matrix. In the iterative procedure, the orthogonal matrix $Y(\tau)$ may be viewed as an improvement of a previous approximation $X$.  

There are many methods available to compute the step-size $\tau$. In \cite[Alg. 2]{Wen13}, the authors recommend a non-monotone linear search method, because of its good theoretical properties regarding the convergence. However, it is not considered here because it has led to poor results in many experiments carried out (not shown here) with particular Moser-Veselov equations. In our modified version of the Wen-Yin algorithm \cite[Alg. 2]{Wen13}, which is presented in Algorithm \ref{alg:wen}, we use, instead, the alternating BB method of \cite{Dai05}: 
\begin{equation}\label{tau}
\tau_k = \left\{\begin{array}{cl}
\|S_{k-1}\|^2_F/|\left< S_{k-1},Y_{k-1}\right>|&\mbox{if}\ k\ \mbox{is odd}\\
|\left< S_{k-1},Y_{k-1}\right>|/\|Y_{k-1}\|^2_F&\mbox{if}\ k\ \mbox{is even}
\end{array}\right.,
\end{equation}
where $S_{k-1}:=X_k-X_{k-1}$, $Y_{k-1}:=\grad F(X_k)-\grad F(X_{k-1})$ and $\left<A,B\right>=\trace(A^TB)$ denotes the Euclidean scalar product. 

\begin{algorithm}[t]
	\caption{Algorithm to solve the Moser-Veselov equation \eqref{eq:moser_veselov} inspired on the steepest descent-type method of Wen\&Yin \cite{Wen13}. The matrix functions $\widetilde{F}(X)$ and $\nabla\widetilde{F}(X)$ are defined, respectively, in (\ref{F-tilde}) and (\ref{grad-F-tilde}).}
	\label{alg:wen}
	{
		\begin{algorithmic}[1]
			\State Choose $\tau>0$ and a starting matrix $X_0\in \mathcal{SO}(n)$;
			\State $f_0\gets \widetilde{F}(X_0)$; $G_0\gets \nabla \widetilde{F}(X_0)$; 
			\State $W_0\gets G_0X_0^T-X_0G_0^T$;
			\State $k\gets 1$;
			\While  {``not converge''}  
			\State $Y_k\gets (I+0.5\tau W_{k-1})^{-1}(I-0.5\tau W_{k-1})X_{k-1}$;
			\State $X_k\gets Y_k$;
			\State $f_{k}\gets \widetilde{F}(X_{k})$ and $G_{k}\gets \nabla \widetilde{F}(X_{k})$; 
			\State $W_{k}\gets G_{k}X_{k}^T-X_{k}G_{k}^T$;
			\State $S_{k}\gets X_k-X_{k-1}$;
			\State $N_k\gets W_k-W_{k-1}$;
			\If {$k$ is even}
			\State $\tau\gets \trace(S_k^TS_k)/|\trace(S_k^TN_k)|$
			\Else
			\State $\tau\gets |\trace(S_k^TN_k)|/\trace(N_k^TN_k)$.
			\EndIf
			\EndWhile
\end{algorithmic}
	}
\end{algorithm}

The second steepest descent-type algorithm addressed here (Algorithm \ref{alg:abrudan}) is a variation of the approaches proposed in \cite{Manton02,Abrudan08} and we refer the reader to those papers for more technical details. In a few words, Algorithm~\ref{alg:abrudan} starts with an initial approximation ${X}_0\in \mathcal{SO}(n)$, finds the skew-symmetric matrix $\text{grad}\, \widetilde F({X_k})$ (the gradient direction on the manifold), and performs several steps along geodesics until convergence. We recall that geodesics on $\mathcal{SO}(n)$ (i.e., curves giving the shortest path between two points in the manifold) can be defined through the matrix exponential as:
$$
G(t)=G(0)\,e^{\mu S},
$$
where ${S}\in\mathbb{R}^{n\times n}$ is a skew-symmetric matrix and $\mu$ is a real scalar.
In Algorithm~\ref{alg:abrudan}, the positive scalar $\mu_k$ controls the length of the ``tangent vector'' and, in turn, the algorithm's overall convergence. To find an almost optimal $\mu_k$, the algorithm uses the Armijo's step-size rule \cite[Sec.1.3]{Polak97}. 

\begin{algorithm}[t]
	\caption{Algorithm to solve the Moser-Veselov equation \eqref{eq:moser_veselov} inspired on the steepest descent-type methods of Manton \cite{Manton02} and Abrudan \cite{Abrudan08}. The matrix functions $\widetilde{F}(X)$ and $\nabla\widetilde{F}(X)$ are defined, respectively, in (\ref{F-tilde}) and (\ref{grad-F-tilde}).}
	\label{alg:abrudan}
	{
		\begin{algorithmic}[1]
			\State ${X}_0\in \mathcal{SO}(n)$ is an initial guess;  
			\State $\mu_1 \gets 1$; 
			\State $\delta \gets 1$; 
			\State $\tau \gets$ {\tt tol};  
			\State $k \gets 0$; 
			\While {$\delta > \tau $} 
			\State ${Z}_k\gets\nabla \widetilde{F}({X}_k)\, {X}_k^T-{X}_k\nabla \widetilde{F}({X}_k)^T$; 
			\State $z_k\gets0.5\ \text{trace}({Z}_k{Z}_k^T)$; 
			\State ${P}_k\gets\text{expm}(-\mu_k{Z}_k)$; 
			\State ${Q}_k\gets{P}_k{P}_k$;  
			\While{$\widetilde{F}({X}_k)-\widetilde{F}({Q}_k{X}_k)\geq \mu_kz_k$}  
			\State ${P}_k\gets{Q}_k$;  
			\State ${Q}_k\gets{P}_k{P}_k$;  
			\State $\mu_k\gets2\mu_k$;  
			\EndWhile
			\While{$\widetilde{F}({X}_k)-\widetilde{F}({Q}_k{X}_k) < 0.5\mu_kz_k$} 
			\State ${P}_k\gets\ \text{expm}(-\mu_k{Z}_k)$;  
			\State $\mu_k\gets0.5 \mu_k$; 
			\EndWhile
			\State ${X}_{k+1}\gets{P}_k{X}_k$; 
			\State $\delta \gets \|{X}_{k+1}-{X}_k\|_{F}$; 
			\State $k\gets k+1$; 
			\EndWhile
			\State ${X}\gets{X}_{k}$. 
		\end{algorithmic}
	}
\end{algorithm}
\section{Numerical Issues}\label{sec:numIssues}
This section addresses some numerical issues associated with the implementation of the three algorithms proposed so far to solve the Moser-Veselov equation.

\subsection{Convergence}

Assuming that all the pure imaginary eigenvalues  (if any) of the $2n\times 2n$ matrix \eqref{eq:hamiltonian} have Jordan blocks with even size, the existence of at least a solution in $\mathcal{SO}(n)$ of the Moser-Veselov equation is guaranteed (check Theorem \ref{thm-existence}). Provided that a careful choice of the starting matrix $X_0$ is made, one of those solutions may be obtained by the three proposed algorithms. We recall that finding an initial guess $X_0$ that minimizes the number of iterations in iterative methods for solving equations is, in general, a challenging problem.  However, steepest-descent algorithms combined with suitable methods for computing the step size have good convergence properties (see \cite{Barzilai88} for the BB method and \cite[Sec. 1.3.2]{Polak97} for the Armijo's method). They have linear convergence, but, in contrast to Newton-type methods, in general, it is easier to find an $X_0$ for which the iterative sequence generated by the method converges. 

In our specific case, we have observed through many experiments (some of them will be shown in Sec.~\ref{experiments}) that $X_0=I$ is a reasonable choice for the three proposed algorithms, in the sense that it leads to convergence towards a special orthogonal solution. Experiments where $X_0$ has been taken as a randomized special orthogonal matrix will be  considered in Sec.~\ref{experiments}, but with a slower convergence; see Figures \ref{fig:experiment2a} and \ref{fig:experiment2b}. 
 
\subsection{Computational Cost}

The three proposed algorithms require $O(n^3)$ operations, which, as written before, is acceptable for algorithms involving matrix computations. However, this information is vague, and we shall give a more sharp estimate of the cost. Namely, we need to find the coefficient of $n^3$ in the polynomial giving the total number of operations. As usual, the terms in $n^2$ and $n$ are ignored.

In terms of computational cost by iteration, the Bregman splitting Algorithm~\ref{alg:bregman} is, in general, the most expensive while Algorithm~\ref{alg:wen} is the cheapest. However, Algorithm~\ref{alg:bregman} converges, in general, faster, attaining the same accuracy of the other two algorithms in much fewer iterations (see Section \ref{experiments}).

\vspace{.25cm}
\noindent
{\bf Algorithm \ref{alg:bregman}:}
The cost is mainly determined by the cost of solving a Sylvester-type equation of the form \eqref{matrix-eq-3} and the computation of two Singular Value Decompositions (SVD). Solving \eqref{matrix-eq-3} involves about $76n^3$ operations ($66n^3$ for the QZ algorithm and $10n^3$ for the remaining calculations; see \cite{Teran11}), and each SVD involves about $22n^3$ operations by the method of Golub and Reinsch \cite{Golub70}. Note that while each iteration includes two SVD, the quite expensive QZ decomposition is just required one time because $A_1$ and $A_2$ are fixed during all the iterations.

\vspace{.25cm}
\noindent
{\bf Algorithms \ref{alg:wen} and \ref{alg:abrudan}:}
The objective functions considered in the steepest descent-type algorithms involve the computation of the trace of matrix products. The efficient computation of $\trace({A}{B})$ does not require matrix-matrix products. Instead, it can be carried out through the formula:
\begin{equation}
\label{eq:trace_simplification}
\trace(AB)= \sum_{i,j}{({A}\circ {B}^T)}_{(i,j)}, 
\end{equation}
where the operator $\circ$ denotes the Hadamard product, i.e., the entry-wise product. If ${A}$ and ${B}$ are matrices of order $n$, the direct computation of the matrix product ${A}{B}$ needs $O(n^3)$ operations, while the trace at \eqref{eq:trace_simplification} just requires $O(n^2)$. However, as far as we know, the trace of a product of three or four matrices requires the computation of one matrix-matrix product which costs $2n^3$ operations. 

Hence, the evaluation of the objective function $\widetilde{F}$  (see its expression in \eqref{F-tilde}) involved in the steepest descent type algorithms requires $4n^3$ operations (two matrix-matrix products), while its computation inside the cycles requires only $2n^3$; the product $JM$ needs to be computed just one time as the algorithm runs. 

Each iteration of the Algorithm~\ref{alg:wen} requires the computation of a Cayley transform, which corresponds to solving a multiple right-hand side linear system of the form $AX=B$, which costs about $8n^3/3$. Concerning Algorithm~\ref{alg:abrudan}, one exponential of a skew-symmetric matrix is required in each iteration. In \cite{Cardoso10}, a scaling and squaring algorithm designed specifically for exponentials of a skew-symmetric matrix is proposed, with an overall cost of $2(16/3+s)n^3$, where $s$ stands for the number of squarings. Alternatively, one can use the general algorithm available through the function \texttt{expm} of MATLAB, which implements the scaling and squaring algorithm of \cite{Mohy09}. Its cost is $O(n^3)$, and we refer the reader to \cite{Mohy09} for the detailed expression of the computational cost. We note that, in the particular case $n=3$, the exponential of a skew-symmetric matrix can be computed by the well-known Rodrigues's formula, at the cost of just one matrix-matrix product.

\subsection{Residual Estimates}

Let us consider the residual function
\begin{equation}\label{res}
R(X):=XJ - JX^T - M
\end{equation}
and assume that $\widetilde{X}$ is an approximation to the exact solution $X$ of the Moser-Veselov equation obtained by a certain numerical algorithm. Hence,  $\widetilde{X}=X+\Delta$, for some matrix $\Delta$ of order $n$.

In the numerical computations of solutions of matrix equations, it is, in general, difficult to estimate the absolute error $\|\Delta\|=\|X-\widetilde{X}\|$ or the relative error $\|\Delta\|/\|X\|$, where $\|.\|$ stands for a given subordinate matrix norm. Thus, the authors work instead with relative residuals to check the quality of the approximation $\widetilde{X}$ and the numerical stability. 

An obvious definition for the relative residual of the Moser-Veselov equation would be 
\begin{equation}\label{rel-res-0}
\|R(\widetilde{X})\|/\|X\|.
\end{equation}
However, as pointed out in \cite[Sec. 5]{Guo06} (see also \cite[Problem 7.15]{Higham08}) for the matrix equation $X^p=A$ (i.e., for computing the matrix $p$th root of $A$), this definition may not be appropriate in some situations. This has motivated us to propose a definition for the relative residual in the style of what is suggested in \cite{Guo06,Higham08}.

With respect to the residual function given in \eqref{res}, we have
\begin{eqnarray*}
	R(\widetilde{X}) &=& (X+\Delta)J-J(X+\Delta)^T-M\\
	&=& XJ+\Delta J-JX^T-J\Delta^T-M\\
	&=& \Delta J-J\Delta^T.
\end{eqnarray*}
By vectorization and attending to some properties of the Kronecker product, we have
\begin{eqnarray*}
	\vec\left( R(\widetilde{X})\right) &=& \vec(\Delta J)-\vec(J\Delta^T)\\
	&=& (J\otimes I)\vec(\Delta)-(I\otimes J)\vec(\Delta^T)\\
	&=& \left(J\otimes I-(I\otimes J)\Pi\right)\vec(\Delta),
\end{eqnarray*}
where $\Pi$ is the permutation matrix of order $n^2\times n^2$ (check (\ref{permut})). Denoting $C:=J\otimes I-(I\otimes J)\Pi\in \mathbb{R}^{n^2\times n^2}$, we have
$$\vec\left( R(\widetilde{X})\right) =  C\vec(\Delta).$$
With respect to the spectral norm $\|.\|_2$, we have
$$\left\|\vec\left( R(\widetilde{X})\right)\right\|_2 \leq \|C\|_2 \|\vec(\Delta)\|_2,$$
and, by attending that, for any matrix $A$, $\|\vec(A)\|_2=\|A\|_F$, it follows
\begin{equation}\label{spectral-ine}    
\left\| R(\widetilde{X})\right\|_F \leq \|C\|_2 \|\Delta\|_F.
\end{equation}
Let us suppose that $\|\Delta\|_F\leq \epsilon \|X\|_F$, for a certain small value $\epsilon$. Note that $\|Q\|_F=\sqrt{n}$, for any orthogonal matrix $Q$ of order $n$ because $\|Q\|_F^2=\trace(Q^TQ)=\trace(I)=\sqrt{n}$. 
Then $\|\Delta\|_F\leq \epsilon \sqrt{n}$ and
$$\|R(\widetilde{X})\|_F\leq \epsilon \sqrt{n} \|C\|_2,$$
or, equivalently,
$$\frac{\|R(\widetilde{X})\|_F}{\sqrt{n}\, \|C\|_2} \leq \epsilon. $$
This suggests the following definition for the relative residual:
\begin{equation}\label{rel-res}
\rho(\widetilde{X}):=\frac{\|R(\widetilde{X})\|_F}{\sqrt{n}\, \|C\|_2}. 
\end{equation}
To understand why the relative residual \eqref{rel-res} is more meaningful than \eqref{rel-res-0}, we shall notice that solving $XJ-JX^T=M$ (ignoring the orthogonal constraint on $X$) is equivalent to solve the linear system
$C\vec(X)=\vec(M)$, with $C$ being in general singular, and that the norm of $C$ must be considered in error analysis, as the expression \eqref{rel-res} does.

\subsection{Termination Criteria}

Several strategies are available to decide when terminating an iterative procedure. Attending to the nature of our optimization problem in \eqref{problem1}, the sequence generated by $F(X_k)$, where 
$$F(X)= \left\|XJ - JX^T - M\right\|^2_F,$$
must converge to zero. Hence, we can fix a tolerance $\epsilon$ and then iterate while $F(X_k)\geq \epsilon$. Note that, in Algorithms~\ref{alg:wen} and \ref{alg:abrudan}, $F(X_k)$ is available during the algorithm and does not involve extra cost. 

The typical behavior of methods with linear convergence is to slow down as $X_k$ approaches stationary points. Often, it may be difficult to detect this phenomenon's occurrence, so another condition to stop the cycle must be added. In our implementations of the algorithm, we consider the classical relative difference
$$\frac{\|X_k-X_{k-1}\|_F}{\|X_k\|_F}=\frac{\|X_k-X_{k-1}\|_F}{\sqrt{n}}$$  
as well. 
Because the algorithms are implemented in finite precision environments, a maximal number of iterations must be considered to stop iterating. In summary, we fix tolerances $\epsilon_1$, $\epsilon_2$ and a maximal number of iterations $k_0$ and stop iterating when
\begin{equation}\label{termination}
\frac{\|X_k-X_{k-1}\|_F}{\sqrt{n}}<\epsilon_2\quad \mbox{or}\quad k > k_0.
\end{equation}
\section{Numerical Experiments}\label{experiments}

To evaluate the performance of the proposed algorithms, we have carried out a set of experiments in MATLAB R2021a (with unit roundoff $u\approx 1.1\times 10^{-16}$) in a machine with Core i5 (1.60GHz). To terminate the iteration procedure in the algorithms, we have used the following criteria:
$$\|X_k-X_{k-1}\|_F/\sqrt{n} \leq \mathtt{tol}\quad \mbox{or}\quad k> 1000,$$ 
where \texttt{tol} is a prescribed tolerance. The following terminology is used:
\begin{itemize}
	\item[$\bullet$] \texttt{\#iter}: number of iterations; 
	\item[$\bullet$] \texttt{rel-res}: relative residual defined in (\ref{rel-res});
	\item[$\bullet$] $\widetilde{F}(X)$: value of the objective function defined in (\ref{F-tilde});
	\item[$\bullet$] $\|\grad\widetilde{F}(X)\|_F$: norm of the Riemannian gradient defined in (\ref{R-grad-F-tilde}).
\end{itemize}
We have considered $r=1$ in Algorithm~\ref{alg:bregman} and $\tau=10^{-3}$ in Algorithm~\ref{alg:wen}.
 Most of the Moser-Veselov equations considered in the experiments do not satisfy the condition $M^2/4+J^2>0$ required by direct methods, but the associated Hamiltonian matrix  ${\mathcal H}$ (see \eqref{eq:hamiltonian}) has no pure imaginary eigenvalue. Experiment 4 is devoted to the case where  ${\mathcal H}$ has pure imaginary eigenvalues. 

\subsection{Experiment 1}

\begin{figure}[t]
    \centering
    \includegraphics[width=0.95\textwidth]{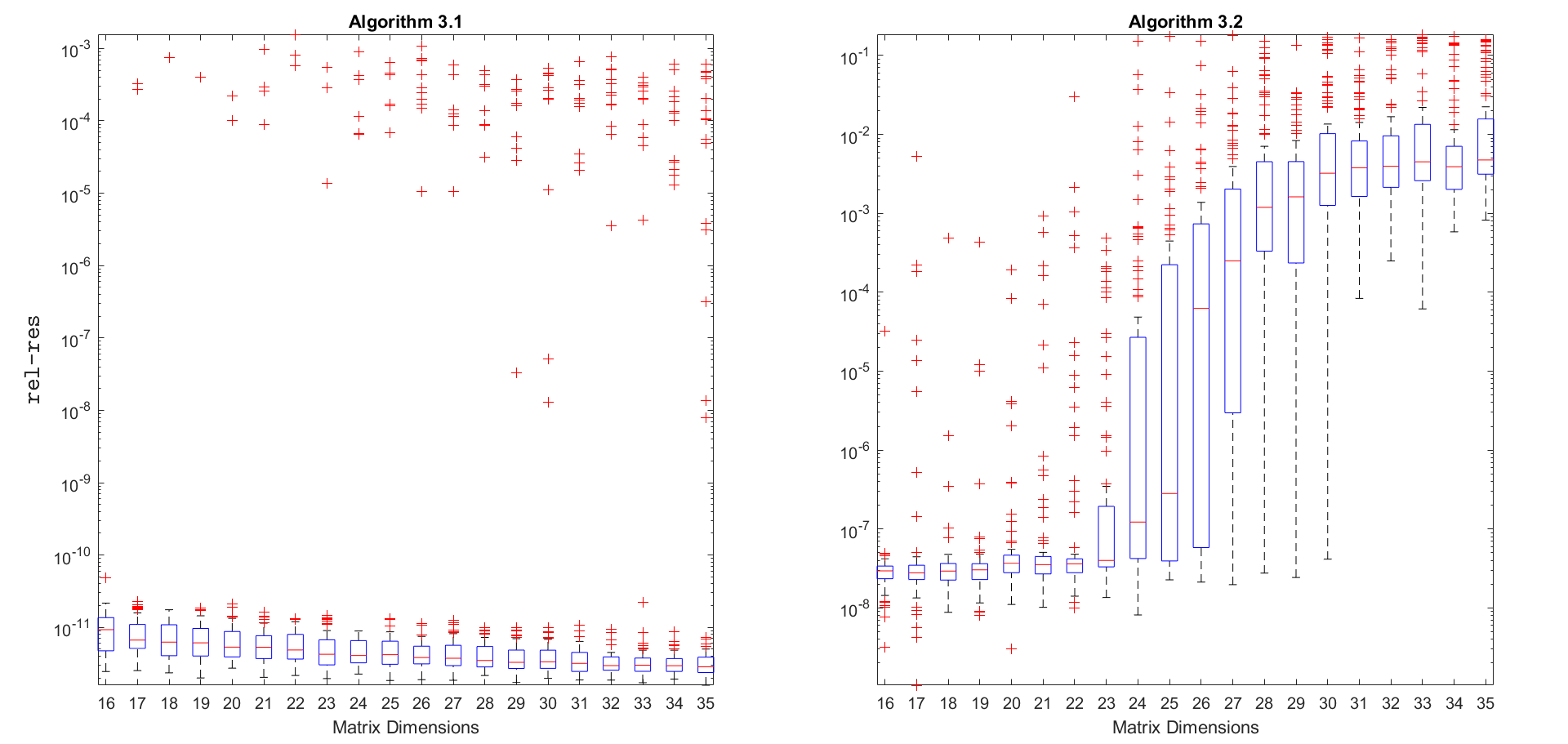}
    \includegraphics[width=0.95\textwidth]{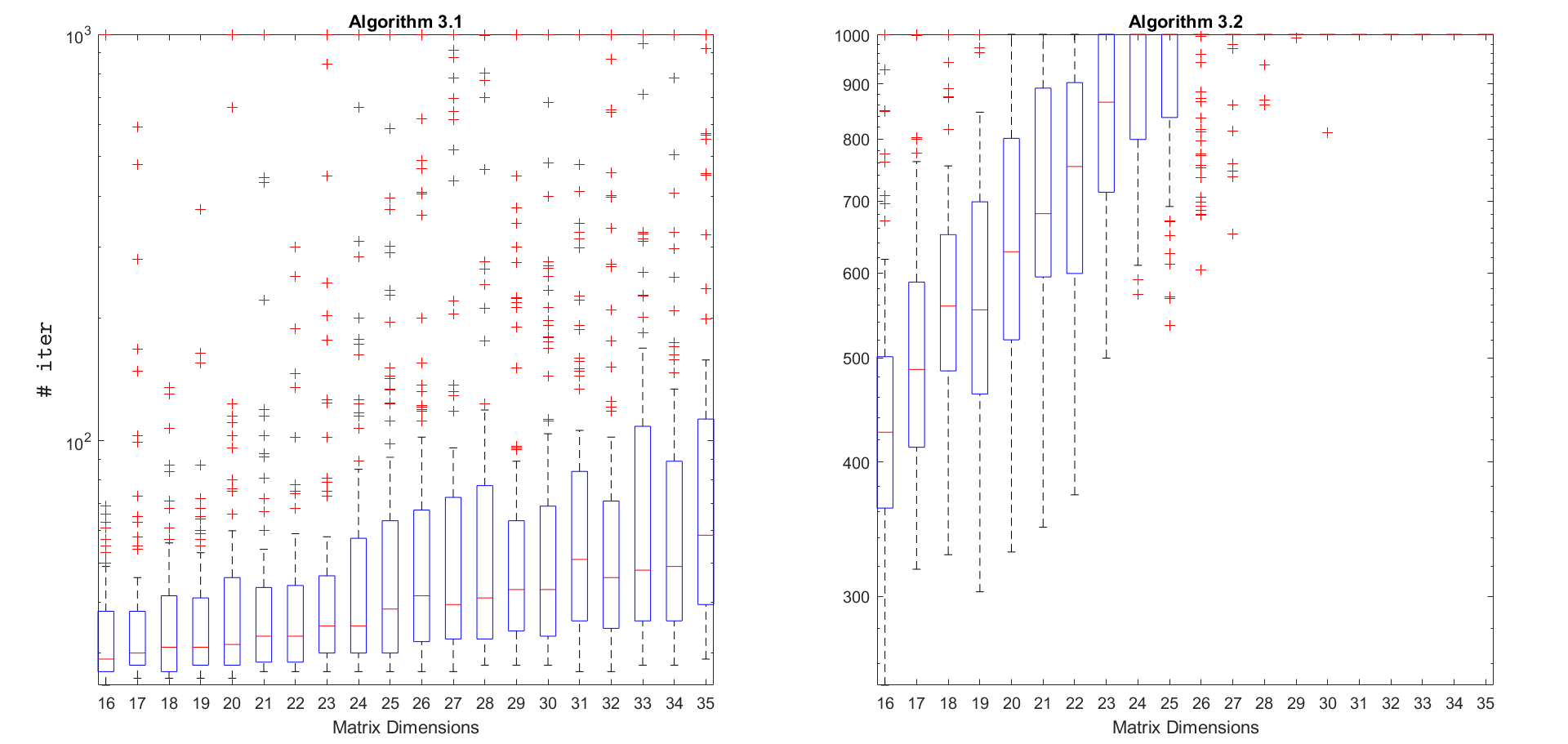}
    \caption{\it Comparison between Algorithm \ref{alg:bregman} and Algorithm \ref{alg:wen} in terms of \texttt{\#iter} and \texttt{rel-res}, for a tolerance $\mathtt{tol}=10^{-10}$ and $X_0=I$.}
    \label{fig:experiment1a}
\end{figure}

\begin{figure}[t] 
    \includegraphics[width=0.95\textwidth]{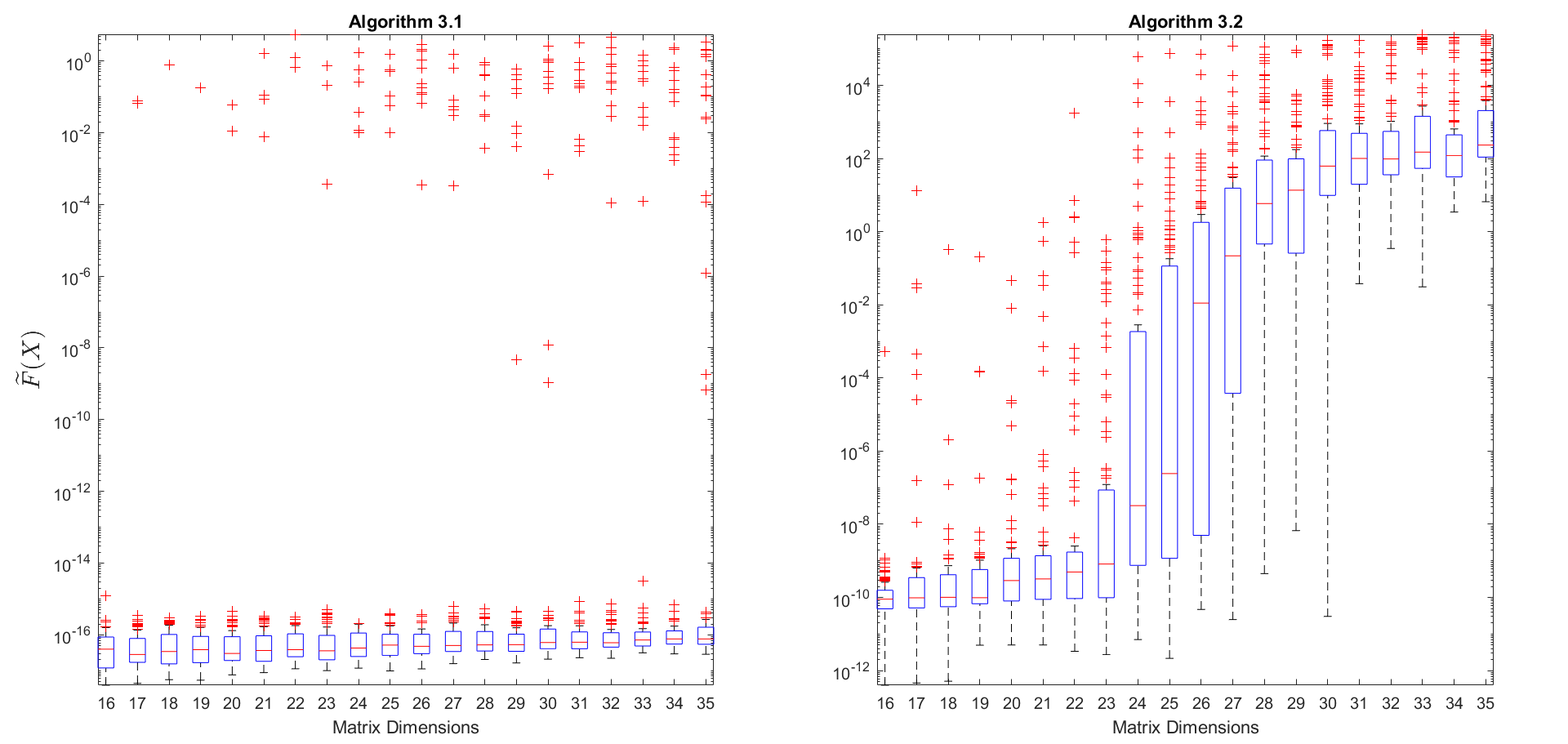}\\
    \includegraphics[width=0.95\textwidth]{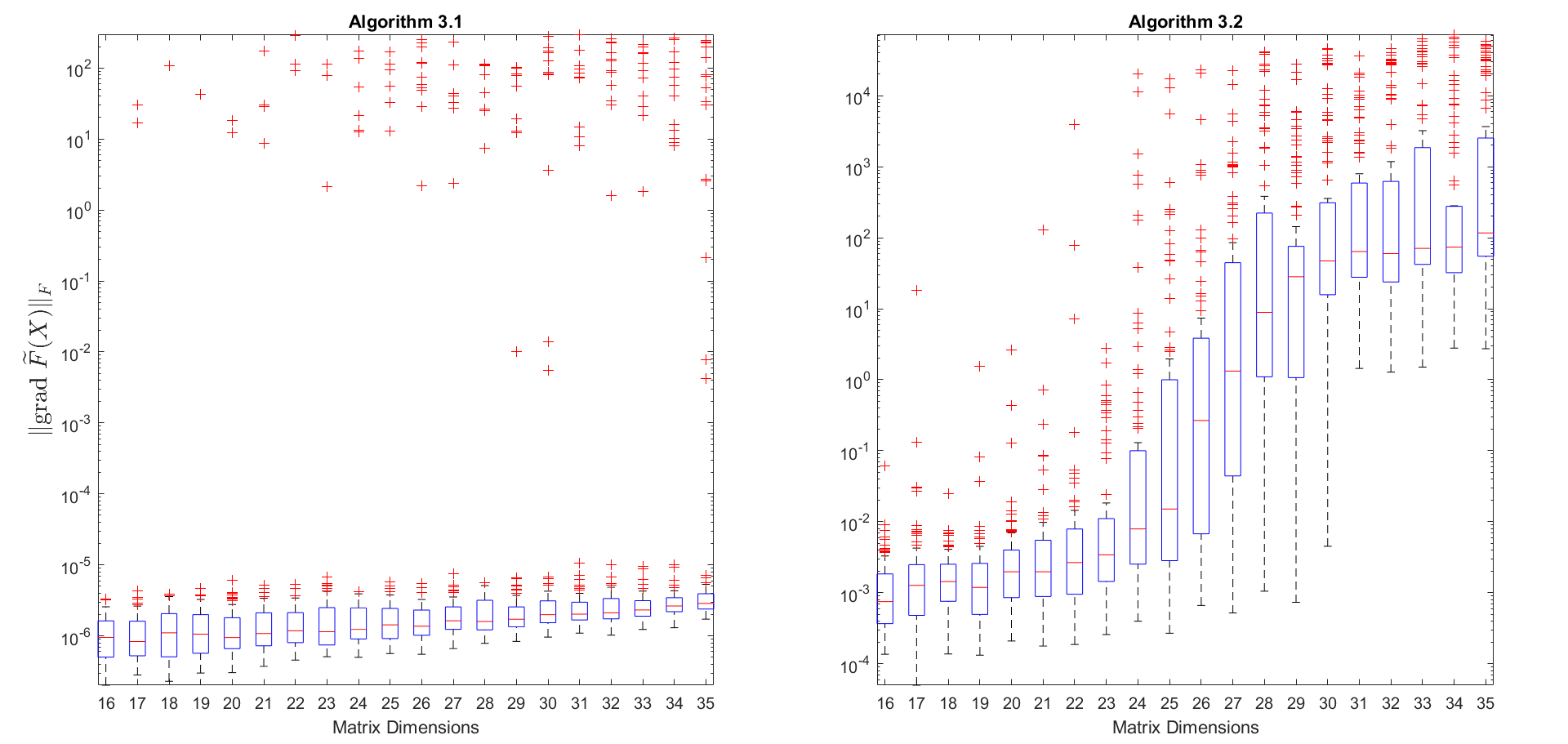}
  \caption{\it Comparison between Algorithm \ref{alg:bregman} and Algorithm \ref{alg:wen} in terms of $\widetilde{F}(X)$ and $\|\grad\widetilde{F}(X)\|_F$, for a tolerance $\mathtt{tol}=10^{-10}$ and $X_0=I$.}
\label{fig:experiment1b}
\end{figure}

This experiment involves a set of $100$ Moser-Veselov equations with randomized matrices $J$ and $M$ of order $16$, plus a set of $100$ Moser-Veselov equations with randomized matrices $J$ and $M$ of order $17$, and so on, up to a set of $100$ Moser-Veselov equations with randomized matrices $J$ and $M$ of order $35$ and aims at comparing Algorithm \ref{alg:bregman} with Algorithm \ref{alg:wen} in terms of \texttt{\#iter}, \texttt{rel-res}, $\widetilde{F}(X)$ and $\|\grad\widetilde{F}(X)\|_F$, for a tolerance $\mathtt{tol}=10^{-10}$ and the initial guess set to $X_0=I$. The results are displayed through a boxplot with whiskers in Figures \ref{fig:experiment1a} and \ref{fig:experiment1b}. 
Experiment 1 involves 2000 Moser-Veselov equations and, as expected, there may have outliers, which are represented by red crosses in the graphs. All of them are finite, i.e., no NaN's of Inf's arose in the calculations. In most of the cases, these outliers reflect difficulties in the intermediate calculations involved in the algorithms, like the computation of inverses of ill-conditioned. 

A careful inspection of those figures lead us to conclude that Algorithm~\ref{alg:bregman}  gives the best results in terms of relative residuals, number of iterations, values of the objective function and norms of the Riemannian gradient, despite it involves higher computational cost by iteration than Algorithm~\ref{alg:wen}. However, since it requires much fewer iterations, the overall computational cost is, in general, smaller. 

\subsection{Experiment 2}\label{sec:experiment2}

\begin{figure}[t]
    \includegraphics[width=0.95\textwidth]{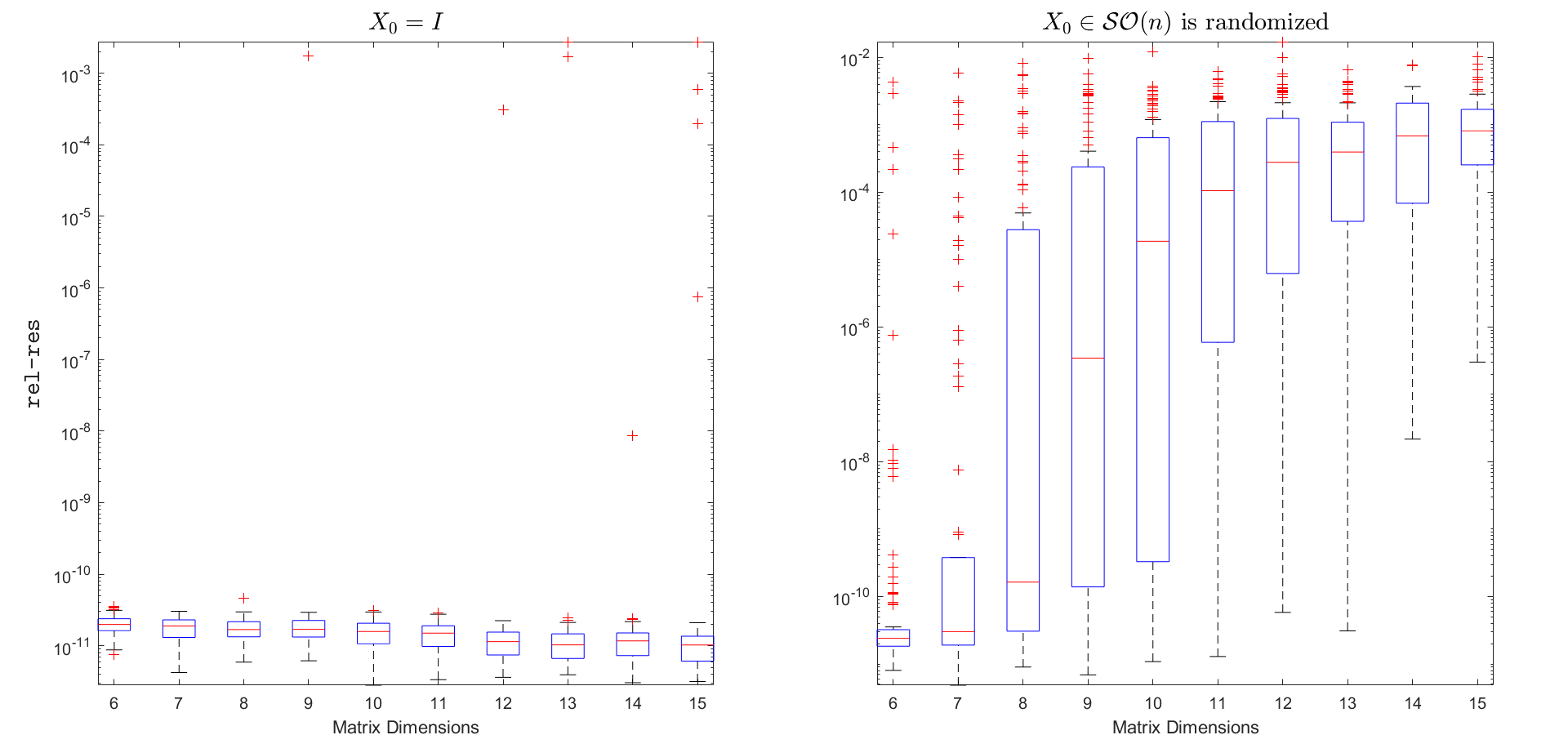}\\
    \includegraphics[width=0.95\textwidth]{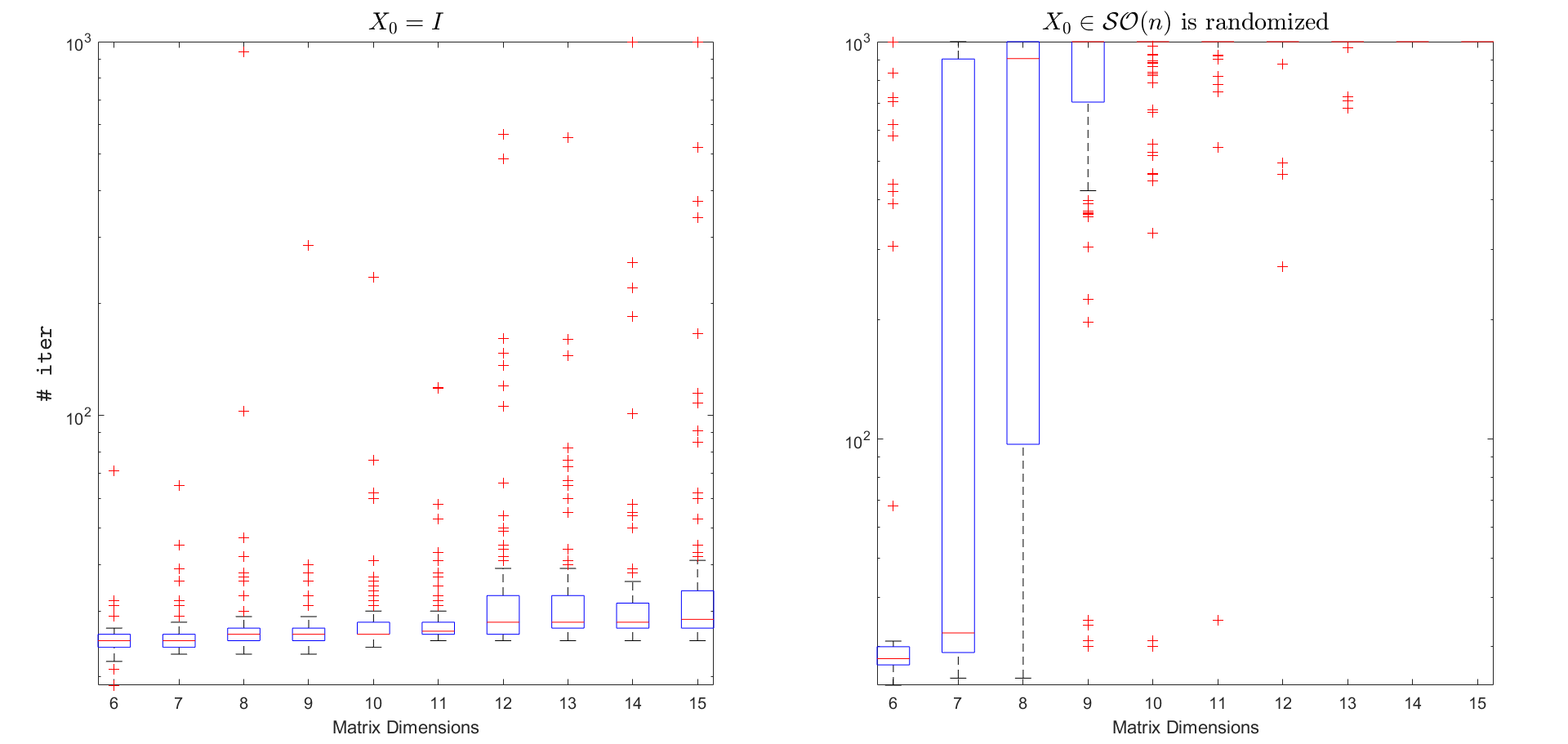}
  \caption{\it Relative residual and number of iterations of Algorithm \ref{alg:bregman} according to the choice of the starting approximation $X_0$ for a tolerance $\mathtt{tol}=10^{-10}$.}
  \label{fig:experiment2a}
\end{figure}

\begin{figure}[t]
    \includegraphics[width=0.95\textwidth]{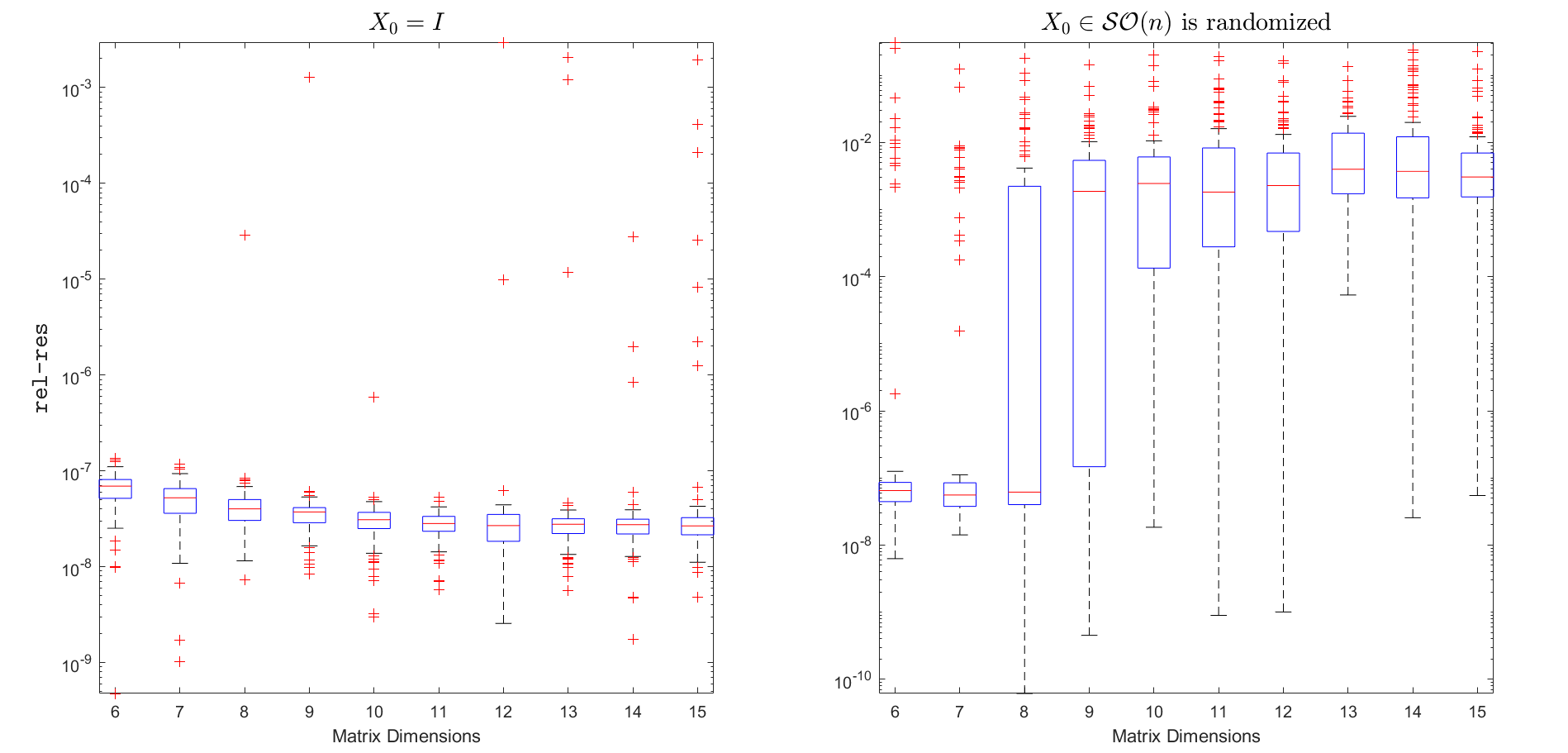}\\
    \includegraphics[width=0.95\textwidth]{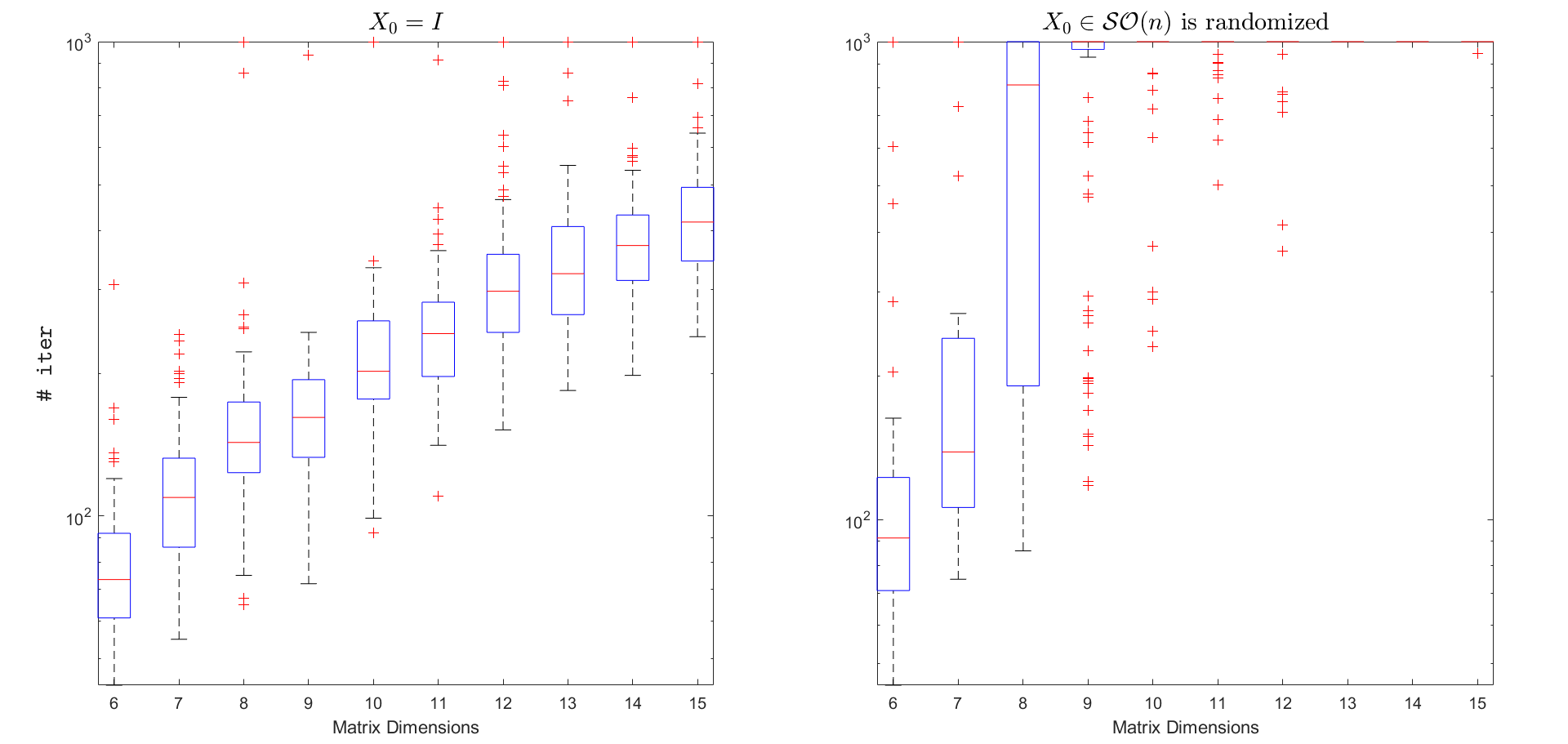}
  \caption{\it Relative residual and number of iterations of Algorithm \ref{alg:wen} according to the choice of the starting approximation $X_0$ for a tolerance $\mathtt{tol}=10^{-10}$.}
  \label{fig:experiment2b}
\end{figure}

This experiment involves $1000$ Moser-Veselov equations: $100$ equations with randomized matrices $J$ and $M$ of order $6$, plus $100$ equations with randomized matrices $J$ and $M$ of order $7$, and so on, up to $100$ equations with randomized matrices $J$ and $M$ of order $15$. It  illustrates the performance of Algorithms \ref{alg:bregman} and \ref{alg:wen} according to the choice of the initial guess $X_0$, in terms of the number of iterations \texttt{\#iter} and of relative residuals \texttt{rel-res}.  The results of Algorithm \ref{alg:bregman} are displayed in Figure \ref{fig:experiment2a}  and of Algorithm  \ref{alg:wen} in Figure \ref{fig:experiment2b},  by means of boxplot graphs with whiskers. We observe that both algorithms perform much better if we take $X_0=I$ instead of a randomized special orthogonal matrix.

\subsection{Experiment 3}

\begin{figure}[t]
    \includegraphics[width=0.95\textwidth]{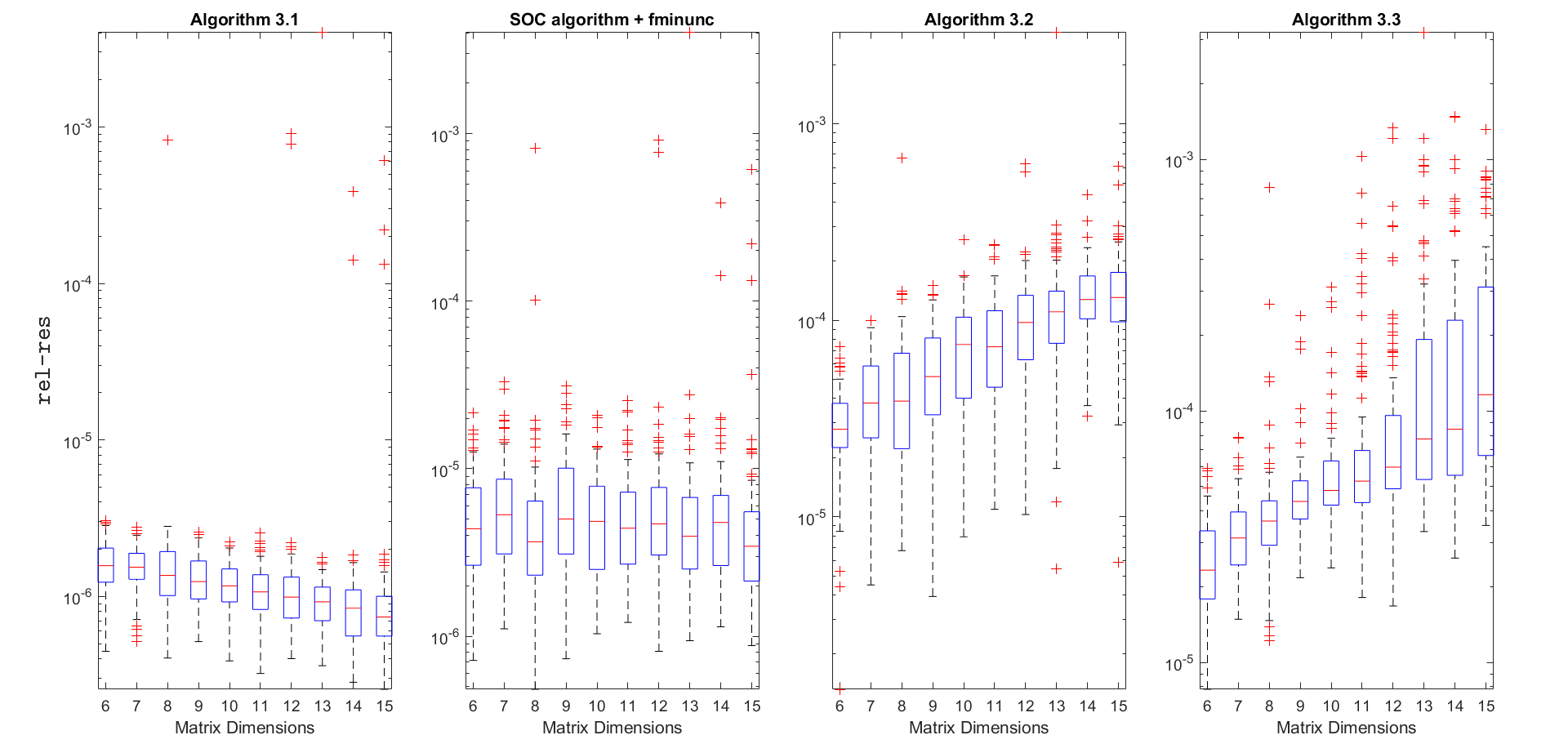}\\
    \includegraphics[width=0.95\textwidth]{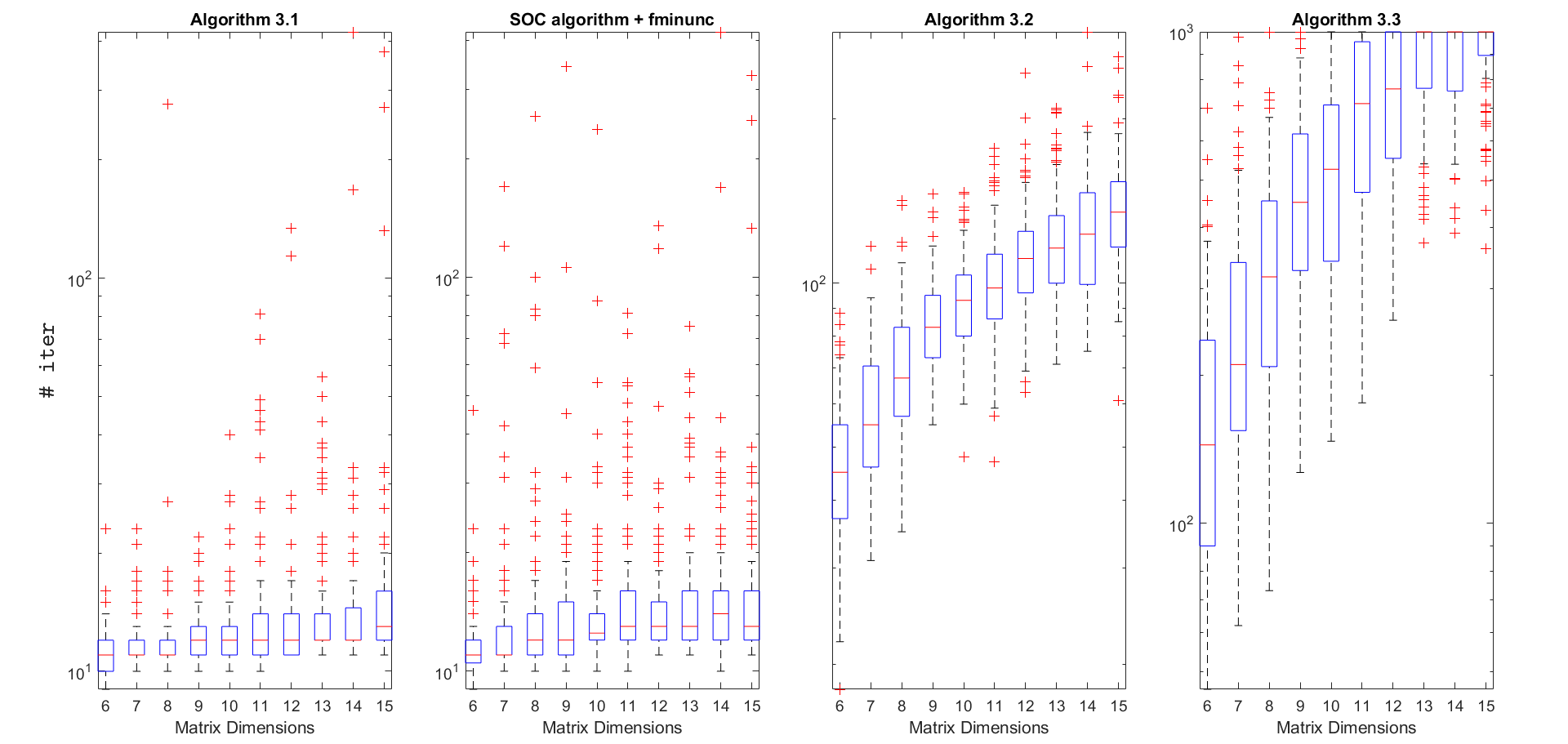}
  \caption{\it Relative residuals and number of iterations of Algorithm \ref{alg:bregman},  SOC Algorithm +\texttt{fminunc}, Algorithm \ref{alg:wen} and Algorithm \ref{alg:abrudan}, for $X_0=I$ and $\mathtt{tol}=10^{-5}$.}
  \label{fig:experiment3a}
\end{figure}

Now, we consider the same set of $1000$ Moser-Veselov equations as in Experiment 2 (Sec. \ref{sec:experiment2}) to illustrate Algorithm \ref{alg:bregman}, SOC Algorithm + \texttt{fminunc} (i.e., the algorithm described at the beginning of Sec.~\ref{sec-alg1}, with Step 2(a) solved using \texttt{fminunc} of MATLAB, where the gradient \eqref{grad-cal-F} is specified in the options mode), Algorithm \ref{alg:wen} and Algorithm \ref{alg:abrudan}, in terms of relative residuals and the number of iterations, for $X_0=I$ and a tolerance $\mathtt{tol}=10^{-5}$. The results are displayed in Figure \ref{fig:experiment3a}. Note that now we are using a larger tolerance than in Experiments 1 and 2. This is because for smaller tolerances both SOC Algorithm+\texttt{fminunc} and Algorithm \ref{alg:abrudan} diverge frequently or may require thousands of iterations. 

\begin{figure}[t]
 \begin{center}
    \includegraphics[width=0.95\textwidth]{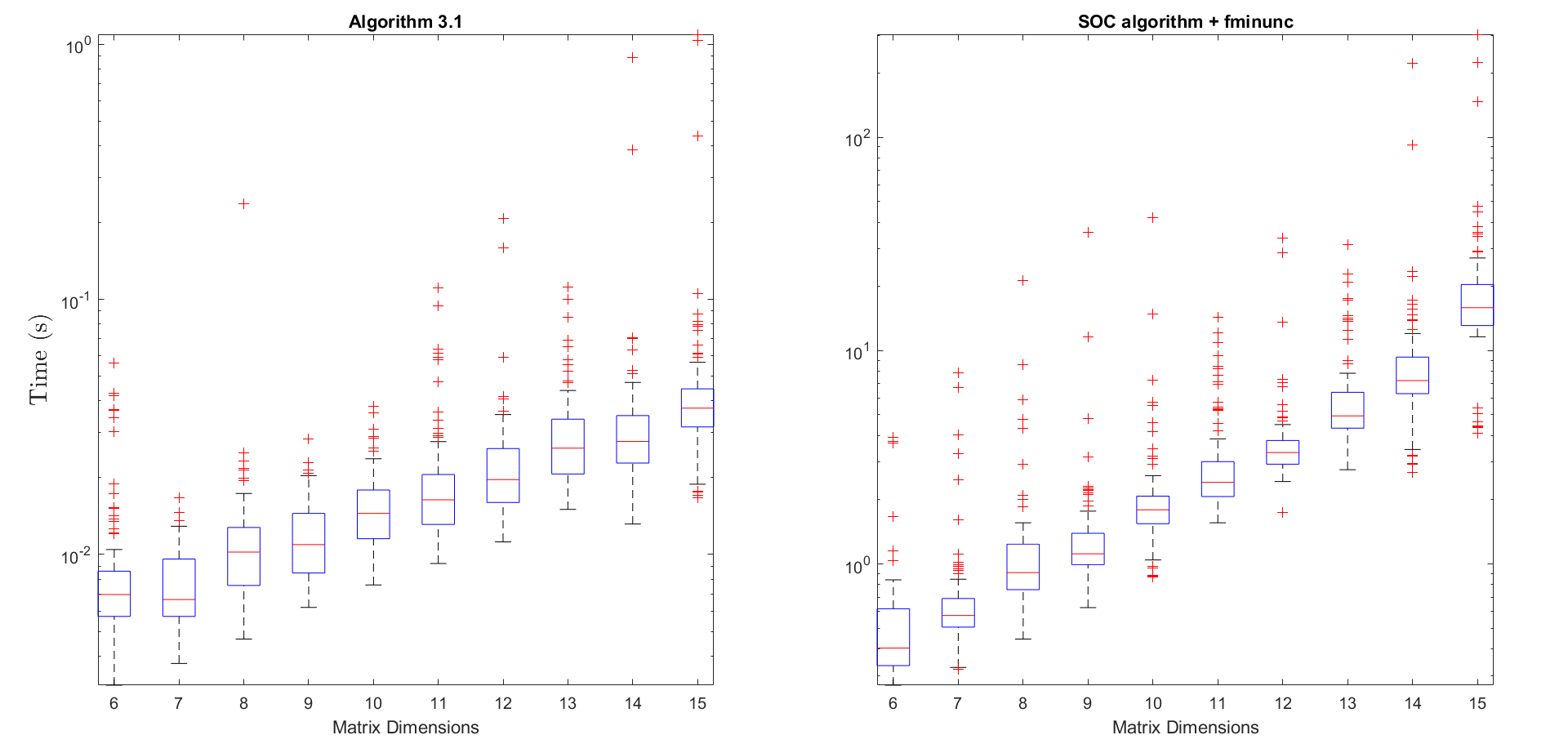}
\end{center}
   \caption{\it Comparison between the computational time of Algorithms \ref{alg:wen} and SOC+\texttt{fminunc}, for $X_0=I$ and $\mathtt{tol}=10^{-5}$.}
  \label{fig:experiment3b}
\end{figure}

In Figure \ref{fig:experiment3b}, we compare the computational time of Algorithm \ref{alg:bregman} with SOC Algorithm + \texttt{fminunc}, where we can see that the latter algorithm requires about $100$ times the computational time of the former. If we decrease the tolerance or increase the size of the matrices, the results of SOC Algorithm + \texttt{fminunc} get even worse and may have little practical interest. For these cases, we recommend instead Algorithms \ref{alg:bregman} or \ref{alg:wen}. 
  
\subsection{Experiment 4}

\begin{figure}[t]
     \includegraphics[width=0.95\textwidth]{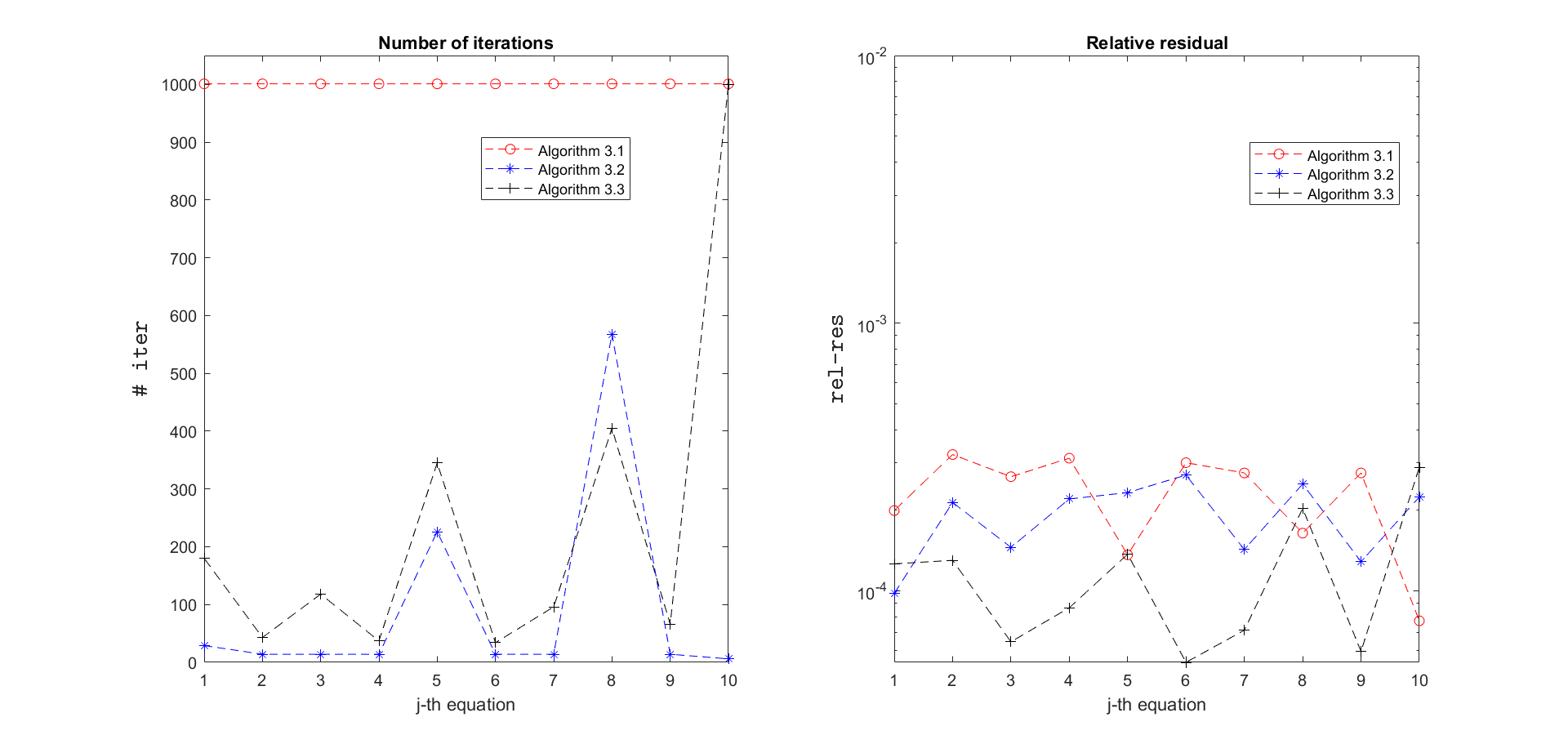}
  \caption{\it Results for Algorithms \ref{alg:bregman}, \ref{alg:wen} and \ref{alg:abrudan} when ${\mathcal H}$ in \eqref{eq:hamiltonian} has some pure imaginary eigenvalues associated to Jordan blocks of even sizes; $X_0=I$ and $\mathtt{tol}=10^{-5}$.}
  \label{fig:experiment4}
\end{figure}

In this experiment, the goal is to illustrate the behaviour of Algorithms \ref{alg:bregman}, \ref{alg:wen} and \ref{alg:abrudan} when the Hamiltonian matrix ${\mathcal H}$ (check \eqref{eq:hamiltonian}) has some pure imaginary eigenvalues associated to Jordan blocks of even sizes.
We have taken ten Moser-Veselov equations involving randomized matrices of order $4$ that were carefully chosen to guarantee that ${\mathcal H}$ fits the conditions just mentioned. The results are displayed in Figure \ref{fig:experiment4}. While in  Experiments 1, 2 and 3, Algorithm \ref{alg:bregman} has performed very well in terms of the number of iterations, now it is the one that gives the poorest results ($1000$ iterations in all the ten tests)! For a compromise regarding the relative residual and the number of iterations, Algorithm \ref{alg:wen} gives the best results. 

\subsection{Other experiments and considerations}

Experiments with Algorithms \ref{alg:bregman} and \ref{alg:wen} for Moser-Veselov equations involving randomized dense matrices of larger sizes ($n=80, 90, 100, 150, 200$) were also performed (not shown here), for a tolerance $\mathtt{tol}=10^{-10}$. In terms of the magnitude of the relative residual and of computational time, we observed a gradual deterioration as the size of the matrices increased. So, our overall recommendation is that those algorithms are suitable for dense matrices of small/medium size, say $n\leq 100$. 

We are not aware of existing algorithms for solving Moser-Veselov equations involving dense (or sparse) matrices of large size.
We leave these investigations for further research.

Methods for solving the optimization problem in \eqref{problem1} hardly give results with relative residuals of order the unit roundoff $u\approx 1.1\times 10^{-16}$. That is the price to pay for considering the objective function as $F(X)=\left\|XJ - JX^T - M\right\|^2_F$ instead of $\left\|XJ - JX^T - M\right\|_F$, whose expression is more difficult to handle. In iterative methods implemented in environments involving floating-point arithmetic, as $F(X)$ becomes less than $u$, its value may stop decreasing. So, it is more reasonable to expect approximations with relative residuals of order $\sqrt{u}$.  

\section{Conclusions} \label{sec:conclusions}

This paper proposes three algorithms for solving the discrete Euler-Arnold equation for the generalized rigid body motion estimation: a Bregman splitting (Algorithm~\ref{alg:bregman}), and two steepest descent-based methods (Algorithms \ref{alg:wen} and \ref{alg:abrudan}). An essential advantage of these methods is that they do not require the strong condition  $M^2/4+J^2>0$, in contrast with other methods existing in the literature. Important numerical issues related to the algorithms, like convergence, computational cost, residual estimates, and termination criteria, have been investigated in detail. The numerical experiments that have been carried out to evaluate the performance of the three methods suggest that the Bregman splitting Algorithm~\ref{alg:bregman} is promising, at least for equations where the associated Hamiltonian matrix has no pure imaginary eigenvalue.

\medskip
\noindent {\bf Acknowledgements:} The work of Jo\~ao R. Cardoso was partially supported by the Centre for Mathematics of the University of Coimbra - UIDB/00324/2020, 
funded by the Portuguese Government through FCT/MCTES. Pedro Miraldo was partially supported by the LARSyS - FCT Plurianual funding 2020-2023.

\section*{References}

\bibliography{refs}

\end{document}